\documentclass[12pt,titlepage]{amsart}
\usepackage{amsfonts, amssymb, verbatim}

\newtheorem*{gentheorem}{Theorem II.1}

\theoremstyle{remark}

\theoremstyle{definition}

\makeatletter
\def\imod#1{\allowbreak\mkern10mu({\operator@font mod}\,\,#1)}
\makeatother
\DeclareMathOperator{\Gal}{{\mathrm{Gal}}}

\DeclareMathOperator{\SL}{{\mathrm{SL}}}

\DeclareMathOperator{\Real}{{\mathrm{Re}}}
 	
\newcommand{\abs}[1]{|#1|}
\newcommand{\elt}[4]{\left(\begin{smallmatrix}#1&#2\\#3&#4\end{smallmatrix}\right)}
\newcommand{\N}{{\mathbf N}}
\newcommand{\Z}{{\mathbf Z}}

\newcommand{\C}{{\mathbf C}}
\newcommand{\HH}{{{\mathbb H}^*}}

\begin{document}
\title{The Divisor Matrix, Dirichlet Series and $\SL(2,\Z)$, II}
\author{Peter Sin}
\address{Department of Mathematics,
University of Florida,
PO Box 118105, Gainesville, FL 32611--8105, USA}
\email{sin@math.ufl.edu}
\author{John G. Thompson}
\address{Department of Mathematics,
University of Florida,
PO Box 118105,
 Gainesville, FL 32611--8105,
 USA}
\email{jthompso@math.ufl.edu}
\date{}
\begin{abstract}
We examine an elliptic curve
constructed in an earlier paper from a certain representation of $\SL(2,\Z)$
on the space of convergent Dirichlet series. The curve is observed to be
a modular curve for $\Gamma^1(15)$ and a certain orbit
of modular functions is thereby associated with the Riemann zeta function. 
Explicit descriptions are given of these functions and
of the permutation action of $\SL(2,\Z)$ on them.
One of the functions has a zero in the arc $\{e^{i\theta}\mid \frac{\pi}{2}\leq\theta\leq\frac{2\pi}{3}\}$. The values of this particular function
along paths formed from $\SL(2,\Z)$-images of the arc 
are used to construct paths between zeros of the zeta function.
\end{abstract}
\maketitle


This paper is a continuation of \cite{ST}. Here we examine some consequences
of Theorems 9.1 and 10.3.

If $X$, $Y$ are indeterminates, set $x=X+1$, $y=Y/X$. We check that
\begin{equation*}
(x-1)y^2+xy-x(x-1)=\frac{Y^2+XY+Y-X^3-X^2}{X}.
\end{equation*}

Since 
\begin{equation*}
(\zeta(s)-1)\phi(s)^2+\zeta(s)\phi(s)-\zeta(s)(\zeta(s)-1)=0, \quad\Real(s)>>0,
\end{equation*}
it follows that
\begin{equation}\label{one}
Y^2+XY+Y=X^3+X^2,
\end{equation}
where $X=\zeta-1$, $Y=\phi/\zeta$.
As the equation (\ref{one}) is in long Weierstrass form \cite[page 2]{SZ},
we compute that it has discriminant $-15$, conductor $15$, and $j$-invariant
$-1/15$. Thus, (\ref{one}) and the Taniyama-Weil Conjecture lead to $\Gamma(15)$.
The curve is labeled 15A in \cite{antwerp} and 15A8 in J. Cremona's tables
(http://www.warwick.ac.uk/staff/J.E.Cremona//ftp/data/).
Following Shimura \cite{Sh} and Fricke \cite{Fr}, if $n\in\N$, set
\begin{equation*}
\Gamma(n)=\lbrace \begin{pmatrix}
a & b\\
c & d\end{pmatrix}\in\SL(2,\Z) \mid a\equiv d\equiv 1, b\equiv c\equiv0\imod n
\rbrace,
\end{equation*}
\begin{equation*}
\Gamma^0(n)=\{\begin{pmatrix}
a & b\\
c & d\end{pmatrix}\in\Gamma(1) \mid  b\equiv0\imod n
\},
\end{equation*}
\begin{equation*}
\Gamma^1(n)=\{\begin{pmatrix}
a & b\\
c & d\end{pmatrix}\in\Gamma^0(n) \mid  a\equiv d\equiv 1\imod n
\}.
\end{equation*}

If $G$ is a subgroup of $\Gamma(1)$ of finite index, $X(G)$ denotes the function
field of $G$.

Following Fricke \cite{Fr}, set
\begin{equation*}
\eta(z)=e^{\pi iz/12}\prod_{n=1}^\infty(1-q^n),\quad q=e^{2\pi iz},\quad \eta_m(z)=\eta(\frac{z}{m}), \quad m\in\N.
\end{equation*}

It is a classical fact that $X(\Gamma^1(15))$ is of genus 1 and has j-invariant $-1/15$.
Therefore, there exist functions $Z$ and $\Phi$ of the upper half plane
which generate $X(\Gamma^1(15))$ as a field extension of $\C$ 
and satisfy the equation 
\begin{equation*}
(Z(z)-1)\Phi(z)^2+Z(z)\Phi(z)-Z(z)(Z(z)-1)=0, \quad \forall z\in\HH.
\end{equation*}

We now proceed to give an explicit construction of such functions $Z$ and $\Phi$.

By Fricke \cite{Fr}, $X(\Gamma^0(15))=\C(\tau,\sigma)$, where
\begin{equation*}
\tau=(\eta_3\eta_5/\eta_1\eta_{15})^3
\end{equation*}
and
\begin{equation}\label{four}
\sigma^2=\tau^4-10\tau^3-13\tau^2+10\tau+1
=(\tau-\alpha)(\tau-\alpha')(\tau-\beta)(\tau-\beta'),
\end{equation}
where
\begin{equation*}
\alpha=\frac{-1}2-\frac{\sqrt{5}}{2},\quad
\alpha'=\frac{-1}2+\frac{\sqrt{5}}{2},\quad
\beta=\frac{11}2-\frac{5\sqrt{5}}{2},\quad
\beta'=\frac{11}2+\frac{5\sqrt{5}}{2}.
\end{equation*}

Since $-I$ fixes every element of $X(G)$ for every subgroup $G$ in $\Gamma(1)$
of finite index, it follows that
$\langle -I, \Gamma^1(15)\rangle$ is the largest subgroup of
$\Gamma(1)$ which acts trivially on $X(\Gamma^1(15))$.
Since $\langle -I, \Gamma^1(15)\rangle\lhd\Gamma^0(15)$ and
$\Gamma^0(15)/\langle -I, \Gamma^1(15)\rangle$ is cyclic of order $4$, it follows that
$X(\Gamma^1(15))$ is a Galois extension of $X(\Gamma^0(15))$ and that
$\Gal(X(\Gamma^1(15))/X(\Gamma^0(15)))\cong\Z/4\Z$.
Let 
\begin{equation*}
\lambda=\eta_1^{-3}\eta_3^6\eta_5^{-3}.
\end{equation*}
Then $\lambda$ is related to $\tau$ and $\sigma$ by the equation
\begin{equation*}
\lambda^2=\frac{c\sigma+d}{250\tau^4},
\end{equation*}
where $c=\sum_{i=0}^4f_i\tau^i$ and $d=\sum_{i=0}^6e_i\tau^i$ with
\begin{multline*}
f_{0}=-1, f_1=-19, f_2=-104, f_3=-125, f_4=125,     \\
e_{0}=1, e_1=24, e_2=180,  e_3=374,  e_4=-396,  e_5=-750, e_6=125.
\end{multline*}

Let $H$ be the subgroup of index 2 in $\Gamma^0(15)$ containing 
$K=\langle -I, \Gamma^1(15)\rangle$. 
Using the standard transformation formulae for $\eta$ 
and generators for $H$, one sees that $\lambda$ is fixed by $H$.
Next, we let $Z$ be a branch of the function 
defined by the equation
\begin{equation}\label{two}
(\tau^2\lambda +\frac{3^35^{-3/2}\tau^3}{\lambda})(Z^2-3Z+1)
=
\sqrt{\frac{5+2\sqrt{5}}{3}}(\tau-\beta)(\tau^2+\gamma\tau+\delta)(Z^2-Z+1),
\end{equation}
with $\gamma=-\frac12-\frac{21}{50}\sqrt{5}$ and $\delta=-\frac{1}{10}-\frac{3}{50}\sqrt{5}$.
To make a definite choice of a germ of $Z$ at $i$, we observe that $\tau(i)$
and $\lambda(i)$ are real, and pick
$Z$ so that $Z(i)$ has positive imaginary part. 
It is straightforward to check that 
\begin{equation}\label{three}
B_1\tau+B_0=0,
\end{equation}
with
\begin{equation*}
B_i=\sum_{j=0}^4 b_{ij}Z^j,\quad i=0,1,
\end{equation*}
\begin{align*}
b_{00}&=b_{04}=1  & b_{01}&=b_{03}=-\frac72-\frac32\sqrt{5},  & b_{02}&=6+3\sqrt{5}\\
b_{10}&=b_{14}=-\frac12-\frac{\sqrt{5}}{2}, & b_{11}&=b_{13}=-2+\sqrt{5}, & b_{12}&=\frac92-\frac32\sqrt{5}.
\end{align*}

Hence,
\begin{equation*}
C(\sigma,\tau, Z)=\C(\sigma,Z).
\end{equation*}
From (2) and (3), we get
\begin{equation*}
B_1^4\sigma^2=(B_0+\alpha B_1)(B_0+\alpha'B_1)(B_0+\beta B_1)(B_0+\beta'B_1).
\end{equation*}
We compute that
\begin{equation*}
B_0+\alpha' B_1= -3\sqrt{5}Z(Z-1)^2,
\end{equation*}
\begin{equation*}
B_0+\alpha B_1= \left(\frac{5+\sqrt{5}}{2}\right)(Z^2-Z+1)^2,
\end{equation*}
\begin{equation*}
B_0+\beta' B_1=-(2+\sqrt{5})(Z-1)^2(4Z^2-7Z+4),
\end{equation*}
\begin{equation*}
B_0+\beta B_1=(\frac92-\frac32\sqrt{5})(Z^2+3Z+1)^2.
\end{equation*}
\begin{equation*}
-3\sqrt{5}\left(\frac{5+\sqrt{5}}{2}\right)(-2-\sqrt{5})(\frac92-\frac32\sqrt{5})=
3^2.5.\left(\frac{1+\sqrt{5}}{2}\right)^2.
\end{equation*}
Set
\begin{equation*}
\Psi=\frac{B_1^2\sigma}{3\sqrt{5}\left(\frac{1+\sqrt{5}}{2}\right)(Z-1)^2(Z^2-Z+1)(Z^2+3Z+1)},
\end{equation*}
so that
\begin{equation}\label{five}
\Psi^2=4Z^3-7Z^2+4Z,
\end{equation}
and
\begin{equation*}
X(\sigma,\tau, Z)=\C(\Psi, Z).
\end{equation*}
Set
\begin{equation}\label{six}
\Phi=\frac{\Psi-Z}{2(Z-1)},
\end{equation}
so that
\begin{equation*}
X(\sigma,\tau, Z)=\C(\Phi, Z).
\end{equation*}
By (\ref{five}) and (\ref{six}),
\begin{equation*}
(Z-1)\Phi^2+Z\Phi-Z(Z-1)=0.
\end{equation*}

We summarize these calculations in the following statement.

\begin{gentheorem}\label{functionfield} 
$X(\Gamma^1(15))$  is  generated as a field extension of $\C$ by the 
functions $Z$ and $\Phi$ defined above, which satisfy 
\begin{equation*}
(Z(z)-1)\Phi(z)^2+Z(z)\Phi(z)-Z(z)(Z(z)-1)=0, \quad \forall z\in\HH.
\end{equation*}
\end{gentheorem}

We compute that
\begin{equation*}
\abs{\Gamma(1):\langle -I, \Gamma^1(15)\rangle}=96,
\end{equation*}
and we construct an explicit set of coset representatives $\{P_i \mid 1\leq i\leq 96\}$,
so that 
\begin{equation*}
\Gamma(1)=\bigcup_{i=1}^{96}\langle -I, \Gamma^1(15)\rangle P_i.
\end{equation*}
These are given in the Appendix.
We set
\begin{equation}\label{seven}
Z_i=Z^{P_i},\quad \mathfrak Z=\{Z_i\mid 1\leq i\leq 96\}.
\end{equation}
The elements of $\mathfrak Z$ are called \emph{avatars} of $\zeta$.

We check that if $z_0\in\HH$ and $Z(z_0)=0$, then
\begin{equation*}
\tau(z_0)=\frac{-1+\sqrt{5}}{2},\quad \sigma(z_0)=0.
\end{equation*}
and from Fricke \cite[page 450, equations (17) and (18)]{Fr},
\begin{equation*}
j=\frac{(\tau_5^2+10\tau_5+5)^3}{\tau_5},\quad \tau_5=(\eta_5/\eta_1)^6,
\end{equation*}

\begin{equation*}
\tau_5=\frac{\tau^4-9\tau^3-9\tau-1+(\tau^2-4\tau-1)\sigma}{2\tau},
\end{equation*}
so
\begin{equation*}
j(z_0)=135\left(\frac{1415+637\sqrt{5}}{2}\right)\approx 629.
\end{equation*}
Set
\begin{equation*}
E=\{e^{i\theta}\mid\frac{\pi}{2}\leq \theta\leq \frac{2\pi}{3}\}.
\end{equation*}
Since $E$ and $[0,1728]$ are in bijection via $z\mapsto j(z)$, $z\in E$,
it follows that there is a unique $c\in E$ such that
\begin{equation*}
j(c)=135\left(\frac{1415+637\sqrt{5}}{2}\right).
\end{equation*}
Since 
$j(z)=j(z')$ if and only if $z$ and $z'$ are in the same $\Gamma(1)$-orbit, it follows
that there is $g_0\in\Gamma(1)$ such that $Z(g_0(c))=0$.
A straightforward calculation yields that
\begin{equation*}
Z^{g_0}=Z_{41}.
\end{equation*}
Set
\begin{equation*}
T=\bigcup_{g\in\Gamma(1)}g(E).
\end{equation*}
We check that if $g\in\Gamma(1)$ and $E\cap g(E)\neq\emptyset$, then one of the following
holds:
\begin{enumerate}
\item[1.] $g=\pm I$ and $E=g(E)$,
\item[2.] $g=\pm R$ and $E\cap g(E)=\{\omega\}$, $\omega=e^{\frac{2\pi i}{3}}$,
\item[3.] $g=\pm S$ and $E\cap g(E)=\{i\}$.
\end{enumerate}
This implies that $T$ is a tree, with vertex set
\begin{equation*}
V(T)=\{g(i) \mid g\in\Gamma(1)\}\cup \{g(\omega) \mid g\in\Gamma(1)\},
\end{equation*}
and edge set
\begin{equation*}
E(T)=\{g(E) \mid g\in\Gamma(1)\}.
\end{equation*}
In particular, if $t_1$, $t_2\in T$, then there is a unique path
from $t_1$ to $t_2$.

It is natural to consider elements $g\in\Gamma(1)$ such that $Z_{41}^g=Z_{41}$
and then to consider
\begin{equation*}
\Delta=\{(z,s)\in\HH\times\C \mid Z_{41}(z)=\zeta(s)\},
\end{equation*}
the motivation being that $Z_{41}$ is an avatar of $\zeta$.

We consider the directed path $P_g$ in $T$ from $c$ to $g(c)$. Thus for each
zero $\rho$ of $\zeta$, $(c,\rho)\in\Delta$.
We then use analytic continuation to build a path $Q_{g,\rho}$ in $\C$
such that as we traverse $P_g$ in $\HH$ we simultaneously traverse $Q_{g,\rho}$
in $\C$ such that if $(z,s)\in\HH\times\C$ and $z$, $s$ are corresponding points, then
$(z,s)\in\Delta$. An injudicious choice of $g$ will have a point $z$ on $P_g$
such that $z$ is a pole of $Z_{41}$, and our way is blocked. However, the element
\begin{equation*}
A=RSRSR^{-1}SRSR^{-1}SRSR^{-1}SRSR=(RSRSR)^4=\begin{pmatrix}0&-1\\1&3\end{pmatrix}^4
\end{equation*}
fixes $Z_{41}$ and $Z_{41}$ has no poles on $P_A$.

Let $\rho_1$, $\rho_2$,\dots be zeros of $\zeta$ in the upper half
plane, ordered so that $\rho_m=\frac12+i\gamma_m$,
with $0<\gamma_1<\gamma_2<\cdots$. Using the computer program SAGE (MAXIMA, Pari) \cite{SAGE}, we found that
$\rho_{m+1}\in Q_{A,\rho_m}$ for $1\leq m\leq 300$. That is,
$(A(c), \rho_{m+1})\in\Delta$, as $Q_{A,\rho_m}$ starts at $\rho_m$
and ends at $\rho_{m+1}$ for $1\leq m\leq 300$. We do not understand this phenomenon, but
it is sufficiently arresting to be noted explicitly.

\section*{Appendix. Coset representatives $P_i$ of $\langle -I, \Gamma^1(15)\rangle$ in $\Gamma(1)$}

{\centering
\begin{tabular}{| c | c | c | c | c |}\hline
$n$&$P_n$&$w_n$&$nR$&$nS$\\[3pt] \hline
$1$ &$\elt1001$ &$1$& $24$ & $24$ \\[3pt]\hline
$2$ &$\elt13{-1}{-2}$ &$R^{-1}SRSR$& $23$ & $20$ \\[3pt]\hline
$3$ &$\elt{-2}31{-2}$ &$R^{-1}SRSR^{-1}S$& $19$ & $23$ \\[3pt]\hline
$4$ &$\elt1{-3}01$ &$R^{-1}SR^{-1}SR^{-1}S$& $44$ & $22$ \\[3pt]\hline
$5$ &$\elt{-2}{-3}11$ &$R^{-1}SR^{-1}SR$& $22$ & $21$ \\[3pt]\hline
$6$ &$\elt{-5}{-3}21$ &$R^{-1}SR^{-1}SRSR^{-1}$& $21$ & $43$ \\[3pt]\hline
$7$ &$\elt{-2}51{-3}$ &$R^{-1}SRSR^{-1}SR^{-1}S$& $36$ & $18$ \\[3pt]\hline
$8$ &$\elt21{-1}0$ &$R^{-1}SR^{-1}$& $11$ & $14$ \\[3pt]\hline
$9$ &$\elt{-4}13{-1}$ &$R^{-1}SRSRSRS$& $55$ & $17$ \\[3pt]\hline
$10$ &$\elt{-1}{-2}11$ &$R^{-1}SR$& $8$ & $16$ \\[3pt]\hline
$11$ &$\elt1{-1}01$ &$R^{-1}S$& $10$ & $15$ \\[3pt]\hline
$12$ &$\elt5{-2}{-2}1$ &$R^{-1}SR^{-1}SRSRS$& $57$ & $13$ \\[3pt]\hline
$13$ &$\elt25{-1}{-2}$ &$R^{-1}SR^{-1}SRSR$& $6$ & $12$ \\[3pt]\hline
$14$ &$\elt1{-2}01$ &$R^{-1}SR^{-1}S$& $5$ & $8$ \\[3pt]\hline
$15$ &$\elt11{-1}0$ &$R^{-1}$& $1$ & $11$ \\[3pt]\hline
$16$ &$\elt{-2}11{-1}$ &$R^{-1}SRS$& $2$ & $10$ \\[3pt]\hline
$17$ &$\elt{-1}{-4}13$ &$R^{-1}SRSRSR$& $28$ & $9$ \\[3pt]\hline
$18$ &$\elt{-5}{-2}31$ &$R^{-1}SRSR^{-1}SR^{-1}$& $3$ & $7$ \\[3pt]\hline
$19$ &$\elt35{-2}{-3}$ &$R^{-1}SRSR^{-1}SR$& $18$ & $30$ \\[3pt]\hline
$20$ &$\elt3{-1}{-2}1$ &$R^{-1}SRSRS$& $17$ & $2$ \\[3pt]\hline
$21$ &$\elt{-3}21{-1}$ &$R^{-1}SR^{-1}SRS$& $13$ & $5$ \\[3pt]\hline
$22$ &$\elt31{-1}0$ &$R^{-1}SR^{-1}SR^{-1}$& $14$ & $4$ \\[3pt]\hline
$23$ &$\elt{-3}{-2}21$ &$R^{-1}SRSR^{-1}$& $16$ & $3$ \\[3pt]\hline
$24$ &$\elt0{-1}11$ &$R$& $15$ & $1$ \\[3pt]\hline
\end{tabular}

\newpage

\begin{tabular}{| c | c | c | c | c |}\hline
$n$&$P_n$&$w_n$&$nR$&$nS$\\[3pt] \hline
$25$ &$\elt{-4}{-15}3{11}$ &$R^{-1}SRSRSR^{-1}SRSRSR$& $48$ & $48$ \\[3pt]\hline 
$26$ &$\elt{-4}31{-1}$ &$R^{-1}SR^{-1}SR^{-1}SRS$& $47$ & $44$\\[3pt] \hline 
$27$ &$\elt{-8}{-3}31$ &$R^{-1}SR^{-1}SRSR^{-1}SR^{-1}$& $43$ & $47$\\[3pt] \hline 
$28$ &$\elt43{-3}{-2}$ &$R^{-1}SRSRSR^{-1}$& $20$ & $46$ \\[3pt] \hline 
$29$ &$\elt73{-5}{-2}$ &$R^{-1}SRSRSR^{-1}SR^{-1}$& $46$ & $45$\\[3pt] \hline 
$30$ &$\elt5{-3}{-3}2$ &$R^{-1}SRSR^{-1}SRS$& $45$ & $19$\\[3pt] \hline 
$31$ &$\elt75{-3}{-2}$ &$R^{-1}SR^{-1}SRSRSR^{-1}$& $12$ & $42$ \\[3pt] \hline 
$32$ &$\elt{-7}45{-3}$ &$R^{-1}SRSRSR^{-1}SRS$& $35$ & $38$  \\[3pt] \hline 
$33$ &$\elt{-7}{-2}41$ &$R^{-1}SRSR^{-1}SR^{-1}SR^{-1}$& $7$ & $65$ \\[3pt] \hline 
$34$ &$\elt{-4}71{-2}$ &$R^{-1}SR^{-1}SR^{-1}SRSR^{-1}S$& $32$ & $40$ \\[3pt] \hline
$35$ &$\elt4{11}{-3}{-8}$ &$R^{-1}SRSRSR^{-1}SRSR$& $34$ & $39$  \\[3pt] \hline 
$36$ &$\elt57{-3}{-4}$ &$R^{-1}SRSR^{-1}SR^{-1}SR$& $33$ & $37$  \\[3pt] \hline 
$37$ &$\elt85{-5}{-3}$ &$R^{-1}SRSR^{-1}SRSR^{-1}$& $30$ & $36$ \\[3pt] \hline 
$38$ &$\elt{-4}{-7}35$ &$R^{-1}SRSRSR^{-1}SR$& $29$ & $32$ \\[3pt] \hline 
$39$ &$\elt{11}{-4}{-8}3$ &$R^{-1}SRSRSR^{-1}SRSRS$& $25$ & $35$ \\[3pt] \hline 
$40$ &$\elt{-7}{-4}21$ &$R^{-1}SR^{-1}SR^{-1}SRSR^{-1}$& $26$ & $34$ \\[3pt] \hline 
$41$ &$\elt41{-1}0$ &$R^{-1}SR^{-1}SR^{-1}SR^{-1}$& $4$ & $81$ \\[3pt] \hline 
$42$ &$\elt58{-2}{-3}$ &$R^{-1}SR^{-1}SRSR^{-1}SR$& $27$ & $31$ \\[3pt] \hline  
$43$ &$\elt{-3}51{-2}$ &$R^{-1}SR^{-1}SRSR^{-1}S$& $42$ & $6$\\[3pt] \hline 
$44$ &$\elt{-3}{-4}11$ &$R^{-1}SR^{-1}SR^{-1}SR$& $41$ & $26$  \\[3pt] \hline 
$45$ &$\elt{-3}{-8}25$ &$R^{-1}SRSR^{-1}SRSR$& $37$ & $29$  \\[3pt] \hline 
$46$ &$\elt3{-4}{-2}3$ &$R^{-1}SRSRSR^{-1}S$& $38$ & $28$  \\[3pt] \hline 
$47$ &$\elt37{-1}{-2}$ &$R^{-1}SR^{-1}SR^{-1}SRSR$& $40$ & $27$  \\[3pt] \hline 
$48$ &$\elt{15}{11}{-11}{-8}$ &$R^{-1}SRSRSR^{-1}SRSRSR^{-1}$& $39$ & $25$  \\[3pt] \hline \end{tabular}
\newpage
\begin{tabular}{| c | c | c | c | c |}\hline
$n$&$P_n$&$w_n$&$nR$&$nS$\\[3pt] \hline
$49$ &$\elt{-7}{-30}4{17}$ &$R^{-1}SRSR^{-1}SR^{-1}SR^{-1}SRSRSRSR$& $72$ & $72$ \\[3pt]\hline 
$50$ &$\elt{-2}91{-5}$ &$R^{-1}SRSR^{-1}SR^{-1}SR^{-1}SR^{-1}S$& $95$ & $68$\\[3pt] \hline 
$51$ &$\elt{-1}{-6}15$ &$R^{-1}SRSRSRSRSR$& $67$ & $71$\\[3pt] \hline 
$52$ &$\elt79{-4}{-5}$ &$R^{-1}SRSR^{-1}SR^{-1}SR^{-1}SR$& $68$ & $70$ \\[3pt] \hline 
$53$ &$\elt{-4}93{-7}$ &$R^{-1}SRSRSRSR^{-1}SR^{-1}S$& $94$ & $69$\\[3pt] \hline 
$54$ &$\elt59{-4}{-7}$ &$R^{-1}SRSRSRSR^{-1}SR$& $69$ & $67$\\[3pt] \hline 
$55$ &$\elt15{-1}{-4}$ &$R^{-1}SRSRSRSR$& $60$ & $66$ \\[3pt] \hline 
$56$ &$\elt{-16}79{-4}$ &$R^{-1}SRSR^{-1}SR^{-1}SR^{-1}SRSRS$& $59$ & $62$  \\[3pt] \hline 
$57$ &$\elt{-2}{-7}13$ &$R^{-1}SR^{-1}SRSRSR$& $31$ & $89$ \\[3pt] \hline 
$58$ &$\elt7{-1}{-6}1$ &$R^{-1}SRSRSRSRSRSRS$& $56$ & $88$ \\[3pt] \hline
$59$ &$\elt7{23}{-4}{-13}$ &$R^{-1}SRSR^{-1}SR^{-1}SR^{-1}SRSRSR$& $58$ & $63$  \\[3pt] \hline 
$60$ &$\elt{-5}{-4}43$ &$R^{-1}SRSRSRSR^{-1}$& $9$ & $61$  \\[3pt] \hline 
$61$ &$\elt{-4}53{-4}$ &$R^{-1}SRSRSRSR^{-1}S$& $54$ & $60$ \\[3pt] \hline 
$62$ &$\elt{-7}{-16}49$ &$R^{-1}SRSR^{-1}SR^{-1}SR^{-1}SRSR$& $77$ & $56$ \\[3pt] \hline 
$63$ &$\elt{23}{-7}{-13}4$ &$R^{-1}SRSR^{-1}SR^{-1}SR^{-1}SRSRSRS$& $49$ & $59$ \\[3pt] \hline 
$64$ &$\elt{-11}{-2}61$ &$R^{-1}SRSR^{-1}SR^{-1}SR^{-1}SR^{-1}SR^{-1}$& $50$ & $82$ \\[3pt] \hline 
$65$ &$\elt{-2}71{-4}$ &$R^{-1}SRSR^{-1}SR^{-1}SR^{-1}S$& $52$ & $33$ \\[3pt] \hline 
$66$ &$\elt5{-1}{-4}1$ &$R^{-1}SRSRSRSRS$& $51$ & $55$ \\[3pt] \hline  
$67$ &$\elt65{-5}{-4}$ &$R^{-1}SRSRSRSRSR^{-1}$& $66$ & $54$\\[3pt] \hline 
$68$ &$\elt{-9}{-2}51$ &$R^{-1}SRSR^{-1}SR^{-1}SR^{-1}SR^{-1}$& $65$ & $50$  \\[3pt] \hline 
$69$ &$\elt{-9}{-4}73$ &$R^{-1}SRSRSRSR^{-1}SR^{-1}$& $61$ & $53$  \\[3pt] \hline 
$70$ &$\elt9{-7}{-5}4$ &$R^{-1}SRSR^{-1}SR^{-1}SR^{-1}SRS$& $62$ & $52$  \\[3pt] \hline 
$71$ &$\elt{-6}15{-1}$ &$R^{-1}SRSRSRSRSRS$& $88$ & $51$  \\[3pt] \hline 
$72$ &$\elt{30}{23}{-17}{-13}$ &$R^{-1}SRSR^{-1}SR^{-1}SR^{-1}SRSRSRSR^{-1}$& $63$ & $49$  \\[3pt] \hline \end{tabular}
\newpage
\begin{tabular}{| c | c | c | c | c |}\hline
$n$&$P_n$&$w_n$&$nR$&$nS$\\[3pt] \hline
$73$ &$\elt{-13}{-30}{10}{23}$ &$R^{-1}SRSRSRSR^{-1}SR^{-1}SR^{-1}SRSR$& $96$ & $96$ \\[3pt]\hline 
$74$ &$\elt{-7}93{-4}$ &$R^{-1}SR^{-1}SRSRSRSR^{-1}S$& $71$ & $92$\\[3pt] \hline 
$75$ &$\elt49{-1}{-2}$ &$R^{-1}SR^{-1}SR^{-1}SR^{-1}SRSR$& $91$ & $95$\\[3pt] \hline 
$76$ &$\elt29{-1}{-4}$ &$R^{-1}SR^{-1}SRSRSRSR$& $92$ & $94$ \\[3pt] \hline 
$77$ &$\elt1{-6}01$ &$R^{-1}SR^{-1}SR^{-1}SR^{-1}SR^{-1}SR^{-1}S$& $70$ & $93$\\[3pt] \hline 
$78$ &$\elt{-5}{-6}11$ &$R^{-1}SR^{-1}SR^{-1}SR^{-1}SR^{-1}SR$& $93$ & $91$\\[3pt] \hline 
$79$ &$\elt{-4}{-5}11$ &$R^{-1}SR^{-1}SR^{-1}SR^{-1}SR$& $84$ & $90$ \\[3pt] \hline 
$80$ &$\elt{-4}{13}3{-10}$ &$R^{-1}SRSRSRSR^{-1}SR^{-1}SR^{-1}S$& $83$ & $86$  \\[3pt] \hline 
$81$ &$\elt1{-4}01$ &$R^{-1}SR^{-1}SR^{-1}SR^{-1}S$& $79$ & $41$ \\[3pt] \hline 
$82$ &$\elt{-2}{11}1{-6}$ &$R^{-1}SRSR^{-1}SR^{-1}SR^{-1}SR^{-1}SR^{-1}S$& $80$ & $64$ \\[3pt] \hline
$83$ &$\elt{13}{17}{-10}{-13}$ &$R^{-1}SRSRSRSR^{-1}SR^{-1}SR^{-1}SR$& $82$ & $87$  \\[3pt] \hline 
$84$ &$\elt51{-1}0$ &$R^{-1}SR^{-1}SR^{-1}SR^{-1}SR^{-1}$& $81$ & $85$  \\[3pt] \hline 
$85$ &$\elt1{-5}01$ &$R^{-1}SR^{-1}SR^{-1}SR^{-1}SR^{-1}S$& $78$ & $84$ \\[3pt] \hline 
$86$ &$\elt{-13}{-4}{10}3$ &$R^{-1}SRSRSRSR^{-1}SR^{-1}SR^{-1}$& $53$ & $80$ \\[3pt] \hline 
$87$ &$\elt{17}{-13}{-13}{10}$ &$R^{-1}SRSRSRSR^{-1}SR^{-1}SR^{-1}SRS$& $73$ & $83$ \\[3pt] \hline 
$88$ &$\elt17{-1}{-6}$ &$R^{-1}SRSRSRSRSRSR$& $74$ & $58$ \\[3pt] \hline 
$89$ &$\elt{-7}23{-1}$ &$R^{-1}SR^{-1}SRSRSRS$& $76$ & $57$ \\[3pt] \hline 
$90$ &$\elt{-5}41{-1}$ &$R^{-1}SR^{-1}SR^{-1}SR^{-1}SRS$& $75$ & $79$ \\[3pt] \hline  
$91$ &$\elt{-9}{-5}21$ &$R^{-1}SR^{-1}SR^{-1}SR^{-1}SRSR^{-1}$& $90$ & $78$\\[3pt] \hline 
$92$ &$\elt{-9}{-7}43$ &$R^{-1}SR^{-1}SRSRSRSR^{-1}$& $89$ & $74$  \\[3pt] \hline 
$93$ &$\elt61{-1}0$ &$R^{-1}SR^{-1}SR^{-1}SR^{-1}SR^{-1}SR^{-1}$& $85$ & $77$  \\[3pt] \hline 
$94$ &$\elt9{-2}{-4}1$ &$R^{-1}SR^{-1}SRSRSRSRS$& $86$ & $76$  \\[3pt] \hline 
$95$ &$\elt9{-4}{-2}1$ &$R^{-1}SR^{-1}SR^{-1}SR^{-1}SRSRS$& $64$ & $75$  \\[3pt] \hline 
$96$ &$\elt{30}{17}{-23}{-13}$ &$R^{-1}SRSRSRSR^{-1}SR^{-1}SR^{-1}SRSR^{-1}$& $87$ & $73$  \\[3pt] \hline \end{tabular}
}


\end{document}